\documentclass[reqno,a4paper,11pt]{article}
 \usepackage{amssymb,amsmath}
 \usepackage{amsfonts}
 \usepackage{mathrsfs,enumerate}
 \usepackage{graphicx}
 \usepackage{diagrams}
 \usepackage{color}
 \usepackage[
 color,
 final] 
 {showkeys}

 \definecolor{refkeybis}{gray}{.65}
 \definecolor{labelkeybis}{gray}{.65}
 {\makeatletter
 \def\SK@refcolor{\color{refkeybis}}%
 \def\SK@labelcolor{\color{labelkeybis}}}

 \numberwithin{equation}{section} 
 \newtheorem{theorem}{Theorem}[section] 
 \newtheorem{lemma}[theorem]{Lemma}
 \newtheorem{corollary}[theorem]{Corollary}
 \newtheorem{proposition}[theorem]{Proposition}
 
 \newtheorem{remark}[theorem]{Remark}

 \newtheorem{definition}[theorem]{Definition}

 \newcommand{\R}{\mathbb{R}}

 \newcommand{\LL}{\mathscr{L}}

 \newcommand{\ee}{{\mbox{\boldmath$e$}}}

 \newcommand{\tauV}{{\kern-3pt\tau}}

 \newcommand{\oVVVk}{\overline{\mbox{\boldmath$V$}}\kern-3pt}
 \newcommand{\tVVVk}{\tilde{\mbox{\boldmath$V$}}\kern-3pt}

 \renewcommand{\restriction}[1]{\lower3pt\hbox{$|_{#1}$}}
 
 \renewcommand{\div}{{\rm div}}
 
 \renewcommand{\to}{\rightarrow}

 \renewcommand{\a}{\alpha}
 \renewcommand{\b}{\beta}

 \newcommand{\e}{\varepsilon}

\newcommand{\dx}{\,{\rm d}x}

\newcommand{\dt}{\,{\rm d}t}
\newcommand{\rd}{{\rm d}}

\def\LL{\mathrm{L}} 

\newcommand{\RR}{\mathbb{R}}

\def\ee{\mathrm{e}} 

\def\qed{\,\unskip\kern 6pt \penalty 500
\raise -2pt\hbox{\vrule \vbox to8pt{\hrule width 6pt
\vfill\hrule}\vrule}\par}
\definecolor{darkblue}{rgb}{0.05, .05, .65}
\definecolor{darkgreen}{rgb}{0.05, .70, .05}
\definecolor{darkred}{rgb}{0.8,0,0}

 \newcommand{\supp}{\operatorname{supp}}
 \newcommand{\sign}{\operatorname{sign}}

 \DeclareMathOperator*{\osc}{osc}
 \DeclareMathOperator*{\esssup}{esssup}
 \DeclareMathOperator*{\essinf}{essinf}

\makeindex
\begin{document}
\title{\bf Total Variation Flow\\and\\
Sign Fast Diffusion\\
in one dimension }
\author{Matteo Bonforte$^{\,a}$
~and~ Alessio Figalli$^{\,b}$}

\maketitle

\begin{abstract}
\noindent We consider the dynamics of the Total Variation Flow (TVF) $u_t=\div(Du/|Du|)$
and of the Sign Fast Diffusion Equation (SFDE)
$u_t=\Delta\sign(u)$ in one spatial dimension.
We find the explicit dynamic and sharp asymptotic behaviour for the TVF,
and we deduce the one for the SFDE by an explicit correspondence between the two equations.
\end{abstract}


%
%

\vfill

\noindent (a) Departamento de Matem\'{a}ticas, Universidad
Aut\'{o}noma de Madrid, Campus de Cantoblanco, 28049 Madrid, Spain.\\
E-mail:\texttt{~matteo.bonforte@uam.es}.
Web-page:\texttt{~http://www.uam.es/matteo.bonforte}\\
\noindent (b) Department of Mathematics, The University of Texas at Austin, 1 University Station C1200,
Austin, TX 78712-1082, USA.\\
E-mail:\texttt{~figalli@math.utexas.edu}.
Web-page:\texttt{~http://www.ma.utexas.edu/users/figalli/}

\section{Introduction}

The Total Variation Flow (TVF) and the Sign-Fast Diffusion Equation (SFDE) are two degenerate
parabolic equations arising respectively as the limit of the (parabolic) $p$-Laplacian as
$p \to 1^+$ and of the fast diffusion equations as $m \to 0^+$.
In one spatial dimension these two equations are strictly related, and our
goal here is to describe their dynamics.

Let us first introduce the SFDE: consider the Fast Diffusion Equation (FDE)
\begin{equation}
\label{eq:FDE}
\partial_t v=\Delta\big(v^m\big),\qquad 0<m<1
\end{equation}
(by definition, $v^m=|v|^{m-1}v$).
By letting $m \to 0^+$ one gets the SFDE
$$
\partial_t v=\Delta\big(\sign(v)\big)\,.
$$

To study the evolution of this equation in one dimension, we will exploit its relation with the TVF
(also called $1$-Laplacian, as it corresponds to the limit of the $p$-Laplacian when $p\to 1^+$): at least formally, if $v$ solves the SFDE, then $u(x):=\int_0^x v$ solves the TVF
\begin{equation}\label{eq:TVF}
\partial_t u= \div\left(\frac{D u}{|D u|}\right).
\end{equation}
 Of course this is purely formal, and we will need to justify it, see Section \ref{sect:TVF fast diff}. The correspondence between solutions of the $p$-Laplacian type equations and solutions to fast diffusion type equations have been used since a long time, cf. \cite{HV}  and more recently in \cite{BaV}, in the study of equations related to the fronts represented in image contour enhancement. In several dimension the correspondence between solutions of the $p$-Laplacian and of the FDE is less explicit and holds only for radial solutions, cf. \cite{ISV}.
 Here our strategy is first to analyze the dynamic of the TVF in one dimension, and then use this to recover the behaviour of solutions of the SFDE.

The literature on the TVF is quite rich, and we suggest to the interested reader the monograph \cite{ACM} as source of references for the existence, uniqueness, basic regularity and different concepts of solutions and their relations, together with estimates on the extinction time (see also the review paper \cite{Giga8} and references therein). The TVF has some interest in applications to noise reduction, cf. \cite{ROF,ACMM,ACM}. The asymptotic behaviour of the TVF is still an open problem in many aspects, even if partial results have appeared in \cite{ACDM, ABC,ACM,BCN,BCN2,BCN3,Giga1}. However,
at least in the simpler case of one spatial dimension, the results contained in the present paper exhibit an almost explicit dynamic and a sharp asymptotic behaviour. After the writing of this paper was essentially completed, we learned of a related work \cite{polacchi}.

\noindent\textbf{Plan of the paper. }

\noindent$\bullet$ In section \ref{TVF1d} we analyze the dynamic of the one-dimensional TVF.
As a first step in this direction, we study the time discretized case,
for which we find the explicit dynamics for ``local step functions'', see Subsection \ref{sols.def.1d} and \ref{sec.step.discr}.
In Subsection \ref{dyn.step.funct} we pass to the continuous time case,
and we find the explicit evolution under the TVF for a generic ``local step functions''.
Then, in Subsection \ref{sect.sols.prop},
exploiting the stability of the TVF in $L^p$ spaces and arguing by approximation
we prove some basic but important properties of solutions to the TVF,
such as the conservation and contractivity property of the local modulus of continuity (Theorem
\ref{local.cont.modulus}), and the explicit behaviour around maxima and minima.

In the case of nonnegative compactly supported initial data,
we first prove an explicit formula for the loss of mass and extinction time for the associated solution to the TVF (Proposition \ref{ext.time+mass.loss}).
Next we study the asymptotic behaviour of such solutions: we analyze and classify the possible asymptotic profiles
(that we shall call more properly ``extinction profiles'', see Theorem \ref{thm.stat.sols}),
and we characterize the asymptotic behaviour of solution to the TVF near the extinction time
as a function of the initial datum (see Theorem \ref{Thm.Asymptotic}).
To our knowledge this is the first completely explicit asymptotic result
for very singular parabolic equations, even if it holds in only one spatial dimension.
Finally, in Subsection \ref{sect.rates} we prove the sharpness of the rate of convergence provided by Theorem \ref{Thm.Asymptotic}.

\noindent$\bullet$ In Section \ref{sect:TVF fast diff} we dedicate our attention to the SFDE.
First we rigorously show that
a BV function solves the TVF if and only if its distributional derivative solves the SFDE,
and then we exploit this describe the dynamic of the SFDE.
We conclude the paper with a discussion on the relation between the SFDE and the Logarithmic Fast Diffusion Equation (LFDE), which is another possible limiting equation of the Fast Diffusion Equation
as $m \to 0^+$, see Subsection \ref{sect:logFDE}.\\

\noindent\textsc{Acknowledgements:} We warmly thank J. L. V\'azquez for useful comments and discussions.
A.F. was partially supported by the NSF grant DMS-0969962. M.B. has been partially funded by Project MTM2008-06326-C02-01  and Ramon y Cajal grant RYC-2008-03521 (Spain).

\section{The 1-dimensional Total Variation Flow}\label{TVF1d}

\normalsize
In this part we deal with the one dimensional TVF.
Before introducing the problem, we first recall briefly some notation and basic facts about BV functions for convenience of the reader.

\subsection{Notations and basic facts about BV functions in one dimension}\label{Notations.1D}
Here we recall some basic facts about one-dimensional BV functions, referring to \cite[Section 3.2]{AFP} for more details.

Consider an open connected interval $I\subseteq\RR$.
A function $u \in BV(I)$ if $u \in L^1_{loc}(I)$, its distributional derivative $Du$ is a (signed) measure, and its total variation $|Du|$
has finite mass.

The distributional derivative $Du$ can be decomposed as $Du=\partial_x u \dx+D^s u$,
where $\partial_x u = \nabla u$ is the absolutely continuous part of $Du$ (with respect to the Lebesgue measure),
and $D^su$ is the singular part.

By Sobolev inequalities we have the inclusion $BV(\RR)\subset L^\infty(\RR)$, and
$BV(I)\cap L^1(I)\subset L^p(I)$ for any $1\leq p \leq \infty$.

If $u \in BV(I)$, up to redefining the function in a set of measure zero,
for every point $x\in I$ it always exists the left or right limit of $u$ at a point $x$, which we denote by
\[
u(x^\pm):=\lim_{y\to x^{\pm}}u(y).
\]
Moreover, the limits above are equal up to a countable number of points.
We will always assume to work with a ``good representative'', so that the above property always holds (see \cite[Theorem 3.28]{AFP}).

\subsection{The setting}\label{sect.sol.TVF}
Let us briefly recall the definition of strong solution to the TVF, in the form we will use it throughout this paper.
For the moment, we do not specify any boundary condition
(so the following discussion could be applied to the Cauchy problem in $\R$, as well as the Dirichlet or the Neumann problem on an interval).

A function $u \in L^\infty([0,\infty),BV(I)) \cap W^{1,2}_{loc}([0,\infty),L^2(I))$ is a strong solution of the TVF if there exists
$z \in L^2_{loc}([0,\infty),W^{1,2}(I))$, with $\|z\|_\infty \leq 1$, such that
\begin{equation}\label{TVF.eq.1d}
\partial_t u= \partial_x z \qquad\mbox{on}\qquad (0,\infty)\times I\,,
\end{equation}
and
\[
\int_0^T\int_I z(t,x)\, Du(t,x)\,dt\,\dx =\int_0^T\int_I |Du(t,x)|\dx\,dt\,\qquad \forall \,T>0.
\]
Roughly speaking, the above condition says that $z=Du/|Du|$.
{We refer to the book \cite{ACM} for a more detailed discussion on the different concepts of solution to the TVF depending on the classes of initial data (entropy solutions, mild solutions, semigroup solution), and equivalence among them.}

{Throughout the paper we will deal with non-negative initial data for the TVF, although many properties maybe extended to signed initial data.}

\subsection{The analysis of the time-discretized problem.}\label{sols.def.1d}
It is well known that the strong solution $u$ of the TVF defined above
is generated via Crandall-Ligget's Theorem, namely it is obtained as the limit of solutions of a time-discretized problem,
formally given by the implicit Euler scheme
\begin{equation*}
\frac{u(t_{i+1})-u(t_i)}{t_{i+1}-t_i}=\partial_x\left(\frac{Du(t_{i+1})}{|Du(t_{i+1})|}\right)\,.
\end{equation*}
{We refer to the book \cite{ACM} fore a more complete and detailed discussion of these facts. }The goal of this section is to understand the behaviour of the time-discretized solution both at continuity and at discontinuity points (see Propositions \ref{lem.cont.points} and \ref{jump.Lemma}).

Let us fix a time step $h>0$, set $t_0=0$, $t_{i+1}=t_i+h=(i+1)h$, and define $u_{ih}(x):=u(ih,x)$
so that $u_0(x)=u(0,x)$. The first step reads:
\begin{equation}\label{h.step}
\frac{u_h(x)-u_0(x)}{h}=\partial_x z_h,
\end{equation}
where $z_h\in L^\infty(I)$ satisfies $\|z_h\|_\infty\le 1\,,\;z_h\,Du_h=\big|Du_h\big|$, and $\partial_x z\in L^2(I)$.
Of course it suffices to understand
the behavior of $u_h$ starting from $u_0$, as all the other steps will follow then by iteration.

We are going to prove the main properties of the time discretized solution,
and to this end is useful to recall an equivalent definition for $u_h$:
\begin{equation}\label{PHI}
u_h={\rm argmin}\big[\Phi_h(u)\big],\qquad\mbox{where}\qquad
    \Phi_h(u)=\int_I|Du| + \frac{1}{2h}\int_I |u-u_0|^2 \dx\,.
\end{equation}
Indeed by strict convexity of the functional $\Phi_h$, the minimizer is unique and is uniquely
characterized by the Euler-Lagrange equation associated to $\Phi_h$, which is exactly \eqref{h.step}, that we shall rewrite in the form
\begin{equation}\label{h.step.1}
u_h= u_0 + h\,\partial_x z_h\,.
\end{equation}
Let us observe that the above construction does not need $u_0$ to be $L^2(I)$:
if $u_0\in BV(I)$ the above scheme still makes sense and provides a function $u_h$ such that $u_h-u_0 \in L^2(I)$.

Next, we remark that since $u_0\,,\,u_h\in BV(I)$ also $\partial_x z_h\in BV(I)\subset L^\infty(I)$,
which implies that $z_h$ is Lipschitz and is differentiable outside a countable set of points. Define the (at most countable) set
\begin{equation}\label{N.z}
N(z_h):=\left\{x\in\RR\;\Big|\; \lim_{\varepsilon\to 0}\frac{z_h(x+\varepsilon)-z_h(x)}{\varepsilon}\;\mbox{does not exists}\right\}\,.
\end{equation}
Since $\partial_x z_h\in BV(I)$, it is continuous outside
$N(z_h)$ (i.e. $\partial_x z_h\in C^0(\RR\setminus N(z_h))$),
and we have that $N(z_h)$ coincides with the set of discontinuity point of $u_h-u_0$\,.

\noindent Collecting all the information obtained so far, we can say that equation \eqref{h.step.1} is equivalent to
\begin{equation}\label{h.infos}
\left\{\begin{array}{lll}
h\,\partial_x z_h(x)=u_h(x)-u_0(x) & \mbox{for all }x\in\RR\setminus N(z_h)\\
|z_h(x)|\le 1\,, & \mbox{for all }x\in\RR\\
z_h(x)=\pm 1\,, & \mbox{for }|Du_h|-\mbox{a.e. }\\
\end{array}\right.
\end{equation}
The next lemmata will allow us to show the important fact that, on $\RR\setminus N(z_h)$, $u_h$ is locally constant
whenever different from $u_0$
(see Proposition \ref{lem.cont.points}).
\begin{lemma}
The following holds:
\begin{equation}\label{meas.N}
|Du_h| \big(\{\partial_x z_h\neq 0\}\setminus N(z_h)\big)=|Du_h|  \big(\{u_h\neq u_0\}\setminus N(z_h)\big)=0\,.
\end{equation}
\end{lemma}
\noindent {\sl Proof.~}The equality
\begin{equation}
\label{eq:dz uh u0}
\{\partial_x z_h\neq 0\}\setminus N(z_h)= \{u_h\neq u_0\}\setminus N(z_h)
\end{equation}
easily follows by observing that $u_h-u_0=h\partial_x z_h \in C^0(\RR\setminus N(z_h))$. Moreover, since $z_h\cdot Du_h= |Du_h|$, we have
\[
\int_{\RR}\big(1-|z_h(x)|\big)\rd|Du_h| =0
\]
which in particular implies
\begin{equation}\label{int.1}
\int_{\{u_h\neq u_0\}\setminus N(z_h)}\big(1-|z_h(x)|\big)\rd|Du_h| =0\,.
\end{equation}
Next we notice that
since $z_h$ is differentiable on $\RR\setminus N(z_h)$ and $|z_h|\leq 1$, we have
$\{z_h=\pm 1\}\cap\{\partial_x z_h\neq 0\}\cap \left(\RR\setminus N(z_h)\right)=\emptyset$.
Hence, thanks to \eqref{eq:dz uh u0}
we deduce that $\big(1-|z_h(x)|\big)>0$ on the set $\{u_h\neq u_0\}\setminus N(z_h)$.
Combining this information with \eqref{int.1} we obtain
that $|Du_h| \big(\{u_h\neq u_0\}\setminus N(z_h)\big)=0$, as desired.\qed\
\

We now use the above lemma to analyze the behavior of $u_h$ near continuity points.
\begin{proposition}[Behaviour near continuity points]\label{lem.cont.points}
If $u_h$ is different from $u_0$ at some common continuity point $x$, then it is constant in an open neighborhood of $x$.
\end{proposition}
\noindent {\sl Proof.~}
Let $x\in I$ be a continuity point both for $u_0$ and $u_h$, and assume that  $u_h(x)>u_0(x)$
(the case  $u_h(x)<u_0(x)$ being analogous).
Then by \eqref{h.infos} we deduce that $\partial_x z_h$ is continuous and strictly positive in an open
neighborhood $I(x)$ of $x$, which together with \eqref{meas.N} implies
$$
|Du_h|\bigl(I(x)\bigr)=0.
$$
Hence $u_h$ is constant on $I(x)$.\qed
\

We now show that if $u_0\in BV(I)$ has some discontinuity jump,
then $u_h$ can only have jumps at such points,
and moreover the size of such jumps cannot increase.
\begin{lemma}[Behaviour at discontinuity points]\label{jump.Lemma}
Let $u_0\in BV(I)$. Then, the following inequalities hold for any $x\in I$:
\begin{equation}\label{jumps.ineq}
\begin{array}{llll}
\mbox{if}\qquad & u_h(x^-) \leq u_h(x^+) & \qquad\mbox{then}\qquad & u_0(x^-)\le u_h(x^-)< u_h(x^+)\le u_0(x^+)\\[3mm]
\mbox{if}\qquad & u_h(x^+)\leq  u_h(x^-) & \qquad\mbox{then}\qquad & u_0(x^+)\le u_h(x^+)< u_h(x^-)\le u_0(x^-)\,.\\
\end{array}
\end{equation}
Moreover,
\begin{equation}\label{jumps.zh}
\begin{array}{llll}
u_h(x^-)< u_h(x^+) & \qquad\mbox{implies}\qquad & z_h(x)=1\\[3mm]
u_h(x^-)> u_h(x^+) & \qquad\mbox{implies}\qquad & z_h(x)=-1\,.\\
\end{array}
\end{equation}
\end{lemma}
\noindent {\sl Proof.~}Let  $x\in I$ be a discontinuity point for $u_h$.
Then
\begin{equation}\label{jump.delta}
Du_h(x)=\big(u_h(x^+)-u_h(x^-)\big)\delta_x\,,
\end{equation}
where $\delta_x$ is the Dirac delta at $x\in I$.  We first prove \eqref{jumps.zh}:
recalling that $z_h \cdot Du_h=\big|Du_h\big|$,
if $u_h(x^-)<u_h(x^+)$ then by \eqref{jump.delta} we get $z_h(x)=1$.
Analogously $u_h(x^-)>u_h(x^+)$ implies $z_h(x)=-1$.

Let us now show \eqref{jumps.ineq}: assume first that $u_h(x^-)<u_h(x^+)$. Since $z_h(x)=1$,
$x$ is a maximum
point for $z_h$, thus $\partial_x z_h(x^-)\ge 0$ and $\partial_x z_h(x^+)\le 0$.
Using \eqref{h.infos}, this implies
\[
\frac{u_h(x^+)-u_0(x^+)}{h}=\partial_x z_h(x^+)\le 0 \le \partial_x z_h(x^-)=\frac{u_h(x^-)-u_0(x^-)}\,{h}
\]
which combined with our assumption $u_h(x^-)<u_h(x^+)$ gives
\[
u_0(x^-)\le u_h(x^-)< u_h(x^+)\le u_0(x^+).
\]
The case $u_h(x^-)>u_h(x^+)$ is analogous.\qed

As an immediate corollary we get:
\begin{corollary}[Local continuity]\label{jump.Coroll}
Let $x \in  I$. If $u_0$ is continuous at $x$, then $u_h$ is continuous at $x$.
\end{corollary}

This result shows that (at least at the discrete level) the TVF cannot create new discontinuities,
and that continuous initial data produce continuous solutions (this is actually what we will prove in
Theorem \ref{local.cont.modulus}).
Observe that this is a local property which does not depend on the boundary conditions
(that we have not specified yet).

We conclude this section with the following estimates on the local loss of mass.
\begin{lemma}[Local $L^1-$estimates]\label{L1.est.sec}
The following estimates hold for any interval $(a,b)\subseteq I$:
\begin{equation}\label{L1.estim}
\left|\int_a^b u_h(x)\dx -\int_a^b u_0(x)\dx\right|\le 2h\,.
\end{equation}
\end{lemma}
\noindent {\sl Proof.~} By equation \eqref{h.infos} we have
\begin{equation}\label{zh.ab}
z_h(b)-z_h(a)=\int_a^b \partial_x z_h(x)\dx =\frac{1}{h}\left[\int_a^b u_h(x)\dx -\int_a^b u_0(x)\dx\right].
\end{equation}
Hence \eqref{L1.estim} follows from the bound $\|z_h\|_\infty\leq 1$.\qed

\subsection{The dynamics of local step functions I. The time discretized case}\label{sec.step.discr}

In this section we use the time discretization scheme to study the dynamics for initial data $u_0$
which coincide with a step function on some open interval $I$.

Let us point out that, if $u_0$ is exactly a step function, then one can give an explicit formula for its evolution (see Section \ref{beh.step.h}  and \cite{Giga1}) by simply checking that it satisfies the equation.
However, by studying the ``time-discretized" evolution
(and then letting $h \to 0$), one can see in a much more natural way the
``locality'' in the dynamic of the TVF.
Moreover, our method shows how to deal with functions which do not belong to $L^2(I)$).
Finally, our description give a good insight of the analysis of the discretized
PDE, which may be useful for numerical purposes.

We would like to notice that it is important for the sequel that the time step $h>0$ is sufficiently small with respect to the size of the jumps of the step function we are considering, as otherwise the discretized dynamics becomes more involved, as we shall show with an example at the end of this section.

\subsubsection{Local evolution of a single step. }
\begin{figure}[ht]
\centering
\includegraphics[height=6.57cm, width=14cm]{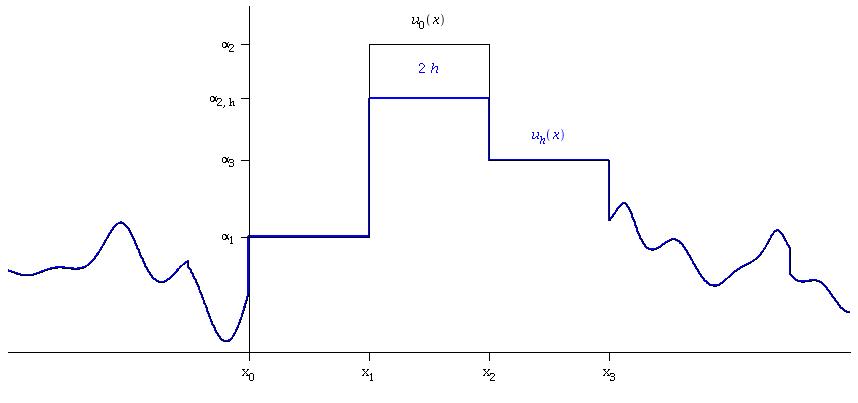}
\caption{\small \textit{Dynamics of a maximum step.
This figure shows the dynamic only inside the interval $[x_0,x_3]$.}}
\label{fig.maxstep}
\end{figure}

To give an insight on the way the discretized evolution behaves, we consider the case of \textit{maximum steps}, whose behavior is made clear in Figure \ref{fig.maxstep}\,. Let us fix an interval $I=I_1\cup I_2\cup I_3$, and assume that $u_0=\alpha_1 \chi_{1} + \alpha_{2}\chi_{2}+ \alpha_{3}\chi_{3}$ on $I$, with $\alpha_2>\max\{\alpha_1,\alpha_3\}$ and $\chi_k=\chi_{I_k}$ is the characteristic function of the open interval $I_k=(x_{k-1},x_k)$.
(Observe that we make no assumptions on $u_0$ outside $I$.)
Fix $h>0$ small (the smallness to be fixed), and consider the function $u_h$.
We have:

\noindent(a) \textit{$u_h$ is constant on any interval $I_k$.}
This follows easily from by Proposition \ref{lem.cont.points} and Corollary \ref{jump.Coroll}.

\noindent(b) \textit{If $h$ is small enough, then $u_h$ jumps at the points $x_1$ and $x_2$.}
Indeed, if by contradiction $x_1$ is a continuity point for $u_h$,
$u_h$ would equal to a constant $\bar \alpha$ on $I_1 \cup I_2$, with
$\alpha_1\le \overline{\alpha}\le \alpha_2$ (see Lemma \ref{jump.Lemma}).
However,
by Lemma \ref{L1.est.sec} we know that
    \[
    |\alpha_1-\overline{\alpha}|\,|I_1|=\left|\int_{I_1} \left[u_h(x)-u_0(x)\right]\dx\right|\le 2h
    \]
    and
    \[
    |\alpha_2-\overline{\alpha}|\,|I_2|=\left|\int_{I_2} \left[u_h(x)-u_0(x)\right]\dx\right|\le 2h
    \]
which is impossible if
    \begin{equation}\label{cond.small.h.orig}
    0<h< \frac{|\alpha_1-\alpha_2|}{2}\left(\frac{1}{|I_1|}+\frac{1}{|I_2|}\right)^{-1}\,.
    \end{equation}
 In particular, $u_h$ jumps both at $x_1$ and $x_2$ if the ``simpler" condition
    \begin{equation}\label{cond.small.h}
    0<h< |\alpha_1-\alpha_2|\,\min\left\{|I_1|\,,\,|I_2|\right\}
    \end{equation}
holds.

\noindent(c) Applying Lemma \ref{jump.Lemma} again (see equation \eqref{jumps.zh}) we obtain
    \[
    z_h(x_1)=1\qquad\mbox{and}\qquad z_h(x_2)=-1\,,
    \]
which combined with \eqref{zh.ab} gives
    \begin{equation*}
    \frac{1}{h}\left[\int_{x_1}^{x_2} u_h(x)\dx -\int_{x_1}^{x_2} u_0(x)\dx\right]=-2\,.
    \end{equation*}
    Hence
    \[
    u_h=\alpha_{2,h}:=\alpha_2-\frac{2h}{|I_2|}\qquad \text{on $I_2$}\,.
    \]

\noindent(d) Combining all together we obtain that
\begin{equation}\label{uh.max.step}
u_0=\alpha_1 \chi_{1} + \alpha_{2}\chi_{2}+ \alpha_{3}\chi_{3}\quad \text{on $I$}\qquad\mbox{implies}\qquad
u_h=\alpha_{1,h} \chi_{1} + \Big(\underbrace{\alpha_{2}-\frac{2h}{|I_2|}}_{\alpha_{2,h}}\Big)\,\chi_{2}+ \alpha_{3,h}\chi_{3}\quad \text{on $I$},
\end{equation}
where $\a_1\leq \a_{1,h} \leq \a_{2,h}$, $\a_3\leq \a_{3,h} \leq \a_{2,h}$
(see Lemma \ref{jump.Lemma}).
The exact value of $\a_{1,h}$ and $\a_{3,h}$ depend on the behavior of $u_0$ outside $I$,
but using  \eqref{zh.ab} we can always estimate them:
\begin{equation}
\label{eq:change steps}
|\a_{1,h}-\a_1| \leq \frac{2h}{|I_1|},\qquad |\a_{3,h}-\a_3| \leq \frac{2h}{|I_3|}.
\end{equation}
(In some explicit cases where one knows that value of $z$ at $x_1$ and $x_3$,
  $\a_{1,h}$ and $\a_{3,h}$ can be explicitly computed using  \eqref{zh.ab}.)

\begin{remark}\label{rmk:local}{\rm
It is important to observe that the value of $\a_{2,h}$ is \emph{independent} of the value of
$\a_{1,h}$ and $\a_{3,h}$, but only depends on the fact that
$\a_1,\a_3 <\a_{2,h}.$ In particular, thanks to \eqref{eq:change steps},
one can iterate the above construction: after $\ell$ steps we get
$$
u_{\ell h}=\alpha_{1,\ell h} \chi_{1} + \alpha_{2,\ell h}\chi_{2}+ \alpha_{3,\ell h}\chi_{3}\quad \text{on $I$},\qquad \a_{2,\ell h}=\alpha_2-\frac{2\ell h}{|I_2|}
$$
holds as long as
$\a_{1,(\ell -1)h},\a_{3,(\ell-1)h} <\a_{2,\ell h},$ which for instance is the case (by iterating the estimate \eqref{eq:change steps}) if
$$
\ell h <|\a_1-\a_2|\min\{|I_1|,|I_2|\} \quad \text{and}\quad \ell h <|\a_2-\a_3|\min\{|I_2|,|I_3|\}.
$$
}
\end{remark}

In particular observe that in this analysis we never used that $I_1$ and $I_3$ are bounded intervals,
so the above formulas also holds when $x_0=-\infty$ and $x_3=+\infty$.

\subsubsection{Evolution of a general step function}\label{beh.step.h}
The above analysis can be easily extended to the general $N$-step function:
assume that
\[
u_0=\sum_{k=0}^{N+1}\alpha_k\chi_k\quad \text{on $I$}
\]
where $\alpha_k\in\RR$ for $k=0,\ldots,N+1$, and $\chi_k=\chi_{I_k}$ is the characteristic function of the open interval $I_k=(x_{k-1},x_k)$ (also the values $x_0=-\infty$ and $x_{N+1}=+\infty$
are allowed). Then, if
\begin{equation}\label{cond.small.h.gen}
0<\ell h<\min_{j=0,\ldots,N}\Big\{\big|\alpha_j-\alpha_{j+1}\big|\min\big\{|I_j|,|I_{j+1}|\big\}\Big\}\,,
\end{equation}
the discrete solution after $\ell$ steps is given by
$$
u_{\ell h}=\sum_{k=0}^{N+1}\alpha_{k,\ell h}\chi_k\quad \text{on $I$},
$$
where we are able to explicitly get the values of $\alpha_{k,\ell h}$ for $k=1,\ldots N$,
see Remark \ref{rmk:local},
and some information on  $\a_{0,\ell k}$ and $\a_{N+1,\ell k}$:
for $k=1,\ldots N$
\begin{equation}\label{alpha.h.k}
\alpha_{k,\ell h}=\left\{
\begin{array}{lll}
\alpha_k\,,\qquad
    &\mbox{if }\alpha_{k-1}< \alpha_k <\alpha_{k+1}  \mbox{ or if }\alpha_{k+1}< \alpha_k <\alpha_{k-1}\\[3mm]
    \alpha_k-\dfrac{2\ell h}{|I_k|}\,,\qquad
    &\mbox{if }\alpha_k>\max\big\{\alpha_{k-1}\,,\,\alpha_{k+1}\big\} \\[3mm]
\alpha_k+\dfrac{2\ell h}{|I_k|}\,,\qquad
    &\mbox{if }\alpha_k<\min\big\{\alpha_{k-1}\,,\,\alpha_{k+1}\big\}\\[3mm]
\end{array}
\right.
\end{equation}
$$
\alpha_{0,\ell h}\left\{
\begin{array}{lll}
\geq \a_{0,(\ell-1)h}\,,\qquad
    &\mbox{if }\alpha_{0}< \alpha_1\\
\leq \a_{0,(\ell-1)h}\,,\qquad
    &\mbox{if }\alpha_{0}> \alpha_1
\end{array}
\right.
$$
$$
\alpha_{N+1,\ell h}\left\{
\begin{array}{lll}
\geq \a_{N+1,(\ell-1)h}\,,\qquad
    &\mbox{if }\alpha_{N}> \alpha_{N+1}\\
\leq \a_{N+1,(\ell-1)h}\,,\qquad
    &\mbox{if }\alpha_{N}< \alpha_{N+1} .
    \end{array}
\right.
$$

\begin{figure}[ht]
\centering
\includegraphics[height=6cm, width=15cm]{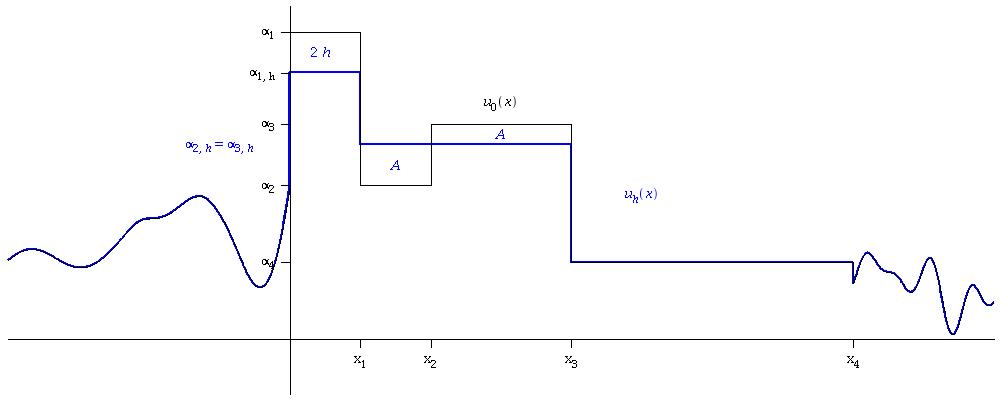}
\caption{\small \textit{Dynamics of a maximum step when $h$ is not necessarily small.
This figure shows the dynamic only inside the interval $[x_0,x_4]$.}}
\label{fig.h.big}
\end{figure}

\noindent\textbf{A concluding remark on the smallness of the time step $h$. }
Since we want to describe the behaviour of the TVF, we are mainly interested in the limit $h\to 0$, which means that condition \eqref{cond.small.h.gen} is always fulfilled. Anyway it is interesting to observe that the dynamic becomes more complicated to understand for general values of $h$,
since the ``locality'' property is lost.
Figure \ref{fig.h.big} shows a situation when a maximum and a minimum disappear in one step
(for this to happen, the area $A$ has to be less than $2h$). Of course one can construct much more complicated examples.
With this one, we can observe that the value of $u_h$ inside
$[x_1,x_2]$ depends on the values of $u_0$ on both $[x_1,x_2]$ and $[x_2,x_3]$.

\subsection{The dynamics of local step functions II. The continuous time case}\label{dyn.step.funct}
We now deduce the explicit evolution of strong solutions to the TVF introduced in Section \ref{sols.def.1d} for the Cauchy problem on $\RR$ where the initial data coincide with a step function (see Remark \ref{rmk:general bdry} below for the analysis of initial value problems with boundary conditions on intervals). The dynamic of step functions will then by used in the next section to deduce, by approximation, qualitative properties for general solutions.

\noindent With the same notation as in Section \ref{beh.step.h}, we consider the initial data given by
\begin{equation}\label{Nstep.initial.data}
u_0(x)=\sum_{k=0}^{N+1}\alpha_k\chi_{I_k}(x) \qquad \text{inside }I.
\end{equation}
Then,  for $h$ small and $x \in I$, we define
\[
u^h(t,x):=\left(\ell+1-\frac{t}{h}\right)u_{\ell h}(x)+\left(\frac{t}{h}-n\right)u_{(\ell+1) h}(x)
    \qquad\mbox{for any }t\in[\ell h, (\ell+1)h]\,.
\]
Then it is immediate to check that, by \eqref{cond.small.h.gen}, on the time interval $[0,t_1]$
with
$$
t_1<\min_{j=0,\ldots,N}\Big\{\big|\alpha_j-\alpha_{j+1}\big|\min\big\{|I_j|,|I_{j+1}|\big\}\Big\}\,,
$$
the solution $u^h$ is explicit, and independent of $h$ on $I_1\cup \ldots\cup I_N$:
$$
u^h(t,x) =u_0(x) + t\,\sum_{k=0}^{N+1} \beta_{k,\ell h}\chi_k(x) \qquad \text{on $[0,t_1]\times I$},
$$
with
\begin{equation}\label{alpha.k.no.h}
\beta_{k,\ell h}:=\left\{
\begin{array}{lll}
0\,,\qquad
    &\mbox{if }\alpha_{k-1}< \alpha_k <\alpha_{k+1}  \mbox{ or if }\alpha_{k+1}< \alpha_k <\alpha_{k-1}\\[5mm]
-\dfrac{2}{|I_k|}\,,\qquad
    &\mbox{if }\alpha_k>\max\big\{\alpha_{k-1}\,,\,\alpha_{k+1}\big\} \\[5mm]
\dfrac{2}{|I_k|}\,,\qquad
    &\mbox{if }\alpha_k<\min\big\{\alpha_{k-1}\,,\,\alpha_{k+1}\big\} \\[5mm]
\end{array}
\right.
\end{equation}
for $k=1,\ldots, N$, and
$$
\beta_{0,\ell h}\left\{
\begin{array}{lll}
\geq 0\,,\qquad
    &\mbox{if }\alpha_{0}< \alpha_1\\
\leq 0\,,\qquad
    &\mbox{if }\alpha_{0}> \alpha_1
\end{array}
\right.
$$
$$
\beta_{N+1,\ell h}\left\{
\begin{array}{lll}
\geq 0\,,\qquad
    &\mbox{if }\alpha_{N}> \alpha_{N+1}\\
\leq 0\,,\qquad
    &\mbox{if }\alpha_{N}< \alpha_{N+1} .
    \end{array}
\right.
$$
By letting $h \to 0$, we find that $u(x,t)$ remains a step function on $I$,
its evolution is explicit on $ \bigl(I_1 \cup \ldots\cup I_{N}\bigr)$,
and on $I_0$ and $I_{N+1}$ it is monotonically increasing/decreasing,
depending on the value on $I_1$ and $I_N$.

This formula will then continue to hold until a maximum/minimum disappear: suppose for instance that $\alpha_2>\max\{\alpha_{1},\alpha_{3}\}$. Then, after a certain time $t'_1$, the value of $u(t)$ on $I_2$ becomes equal to $\max\{\alpha_{1},\alpha_{3}\}$.
Then, we simply take $u(t'_1)$ as initial data and we repeat the construction.

After repeating this at most $N$ times, all the maxima and minima inside $I$ disappear, and $u(t)$
is monotonically decreasing/increasing on $I$.

For instance, if $I=\R$ and the initial data is a compactly supported step function, then
$u\equiv 0$ after some finite time $T$ (which we call \emph{extinction time}).
On the other hand, if $u_0$ is an increasing (resp. decreasing) step function, then it will remain constant in time.

\begin{remark}\label{rmk:general bdry}
{\rm
The analysis done up to now can be extended to the case of suitable initial value problems on intervals with boundary condition. For instance, the dynamic of the Dirichlet problem is analogous to the one described above for the Cauchy problem with compactly supported initial data; we leave the details to the interested reader. Let us consider next the Neumann problem on some closed interval $[a,b]=I_0\cup \ldots
\cup I_{N+1}$. The dynamics on $I_1 \cup \ldots \cup I_N$ is known by our analysis (which, as we observed before, is ``local'').
To understand the dynamics on $I_0$ and $I_N$, we go back to the time discretized problem:
the Euler-Lagrange equations in this case are still \eqref{h.infos}, but with the additional Neumann
condition $z_h(b)=z_h(a)=0$. It is easy to check that this last condition allows to uniquely characterize the value of $u_h$ inside $I_0$ and $I_{N+1}$.

For example, if $u_0=\sum_{k=1}^{N+1} \a_k \chi_{I_k}$ with $\a_1 \leq \ldots\leq \a_{N+1}$
(i.e. $u_0$ is monotonically increasing), then
$$
u(t)=u_0+t \left(\frac{1}{|I_0|}\chi_{I_0}-\frac{1}{|I_{N+1}|}\chi_{N+1}\right)
$$
(i.e. the value on $I_0$ increases, while the one on $I_{N+1}$ decreases).
This holds true until a jump disappears, and then one simply repeat the construction. We leave the details for the general case to the interested reader.
}
\end{remark}


\subsection{Some properties of solutions to the TVF}\label{sect.sols.prop}
 In this section we prove some local ``regularity'' properties enjoyed by solutions of the TVF (see also \cite{CCN, GigaARMA}). The key fact behind these results is that the TVF is contractive in any $L^q$ space with $q \in [1,\infty]$, cf. for example \cite{ACM}. This contractivity property is not so surprising, since it holds also for the $p$-Laplacian for any $p>1$. We are going to show that most of the properties which holds in the case of step functions, can be easily extended to the general case. An example is the following:

\begin{theorem}[Local continuity]\label{local.cont.modulus}
Assume  that $u_0$ is continuous on some open interval $I$. Then also the corresponding solution $u(t)$ is continuous on the same interval $I$ and the oscillation is contractive, namely
\[
\sup_I u(t)-\inf_I u(t)=:\osc_I\big(u(t)\big)\le \osc_I\big(u_0\big).
\]
\end{theorem}
\noindent {\sl Proof.~}
By the contractivity of the TVF in $L^\infty$,  we have
\begin{equation}\label{L.infty.contr}
\|u(t)-v(t)\|_\infty\le \|u_0-v_0\|_\infty
\end{equation}
for any two given solutions $u(t), v(t)$ corresponding to initial data $u_0,v_0$, with $u_0-v_0\in L^\infty$.

Since $u_0$ is continuous in $I$, we can find a family of functions $u_0^\varepsilon$ such that $u_0^\varepsilon=u_0$ outside $I$, $u_0^\varepsilon$ are step functions inside $I$,
and $\|u_0-u_0^\varepsilon\|_\infty\le \varepsilon$. Then by \eqref{L.infty.contr} we get
\[
\|u(t)-u^\varepsilon(t)\|_\infty\le \|u_0-u_0^\varepsilon\|_\infty\le \varepsilon\qquad\mbox{for all }t>0
\]
where $u^\varepsilon(t)$ is the solution to the TVF corresponding to $u_0^\varepsilon$, which is still a step function inside $I$. We now observe that the elementary inequality
\begin{equation}\label{osc.est.1}
\big|\osc_I(u(t))-\osc_I(u^\varepsilon(t))\big|\le 2\,\|u(t)-u^\varepsilon(t)\|_{L^\infty(I)} \leq 2 \varepsilon
\end{equation}
holds. Moreover,
by looking at the explicit formulas for the evolution of $u^\varepsilon(t)$ inside $I$, cf. Section \ref{dyn.step.funct}, it is immediate to check that  $\osc_I(u^\varepsilon(t))$ is decreasing in time.
Hence
$$
\osc_I \bigl(u(t)\bigr) \leq \osc_I\bigl(u^\varepsilon(t)\bigr) +2\varepsilon \leq
\osc_I\bigl(u_0^\varepsilon\bigr) +2\varepsilon \leq
\osc_I\bigl(u_0\bigr) +4\varepsilon.
$$
We conclude letting $\varepsilon\to 0$.\qed

\begin{remark}{\rm
The above theorem still holds if $u_0$ is not continuous on $I$:
in that case one has to replace sup and inf by esssup and essinf, and to prove the result one can
use the comparison principle: if $u^+(t)$ and $u^-(t)$ are the solution starting respectively from
$$
u^+(x):=\left\{
\begin{array}{ll}
u_0(x)&\text{if }x \not\in I;\\
\esssup\limits_I u_0 &\text{if }x \in I;\\
\end{array}
\right.
\qquad
u^-(x):=\left\{
\begin{array}{ll}
u_0(x)&\text{if }x \not\in I;\\
\essinf\limits_I u_0 &\text{if }x \in I;\\
\end{array}
\right.
$$
then $u^-(t) \leq u(t)\leq u^+(t)$, $u^+(t)$ and $u^-(t)$ are both constant on $I$,
and $\|u^-(t,x) - u^-(t,x)\|_{L^\infty(I)}$ is decreasing in time.
However, since we will never use this fact, we leave the details to the interested reader.}
\end{remark}

\subsection{Further properties}\label{sect:further}
Arguing by approximation as done in Theorem \ref{local.cont.modulus} above
(using either the stability in $L^\infty$
or simply the stability in $L^1$, depending on the situation), we can easily deduce other
\textit{local properties of the TVF, valid on any subinterval} $I$ where the solution $u(t)$ is considered (we leave the details of the proof to the interested reader):

\noindent(i) The set of discontinuity points of $u(t)$ is contained in the set of discontinuity points of $u_0$, i.e. ``the TVF does not create new discontinuities''.

\noindent(ii) The number of maxima and minima decreases in time.

\noindent(iii) If $u_0$ is monotone on an interval $I$, then $u(t)$ has the same monotonicity as $u_0$ on $I$. Moreover, if $u_0$ is monotone on the whole $\R$, then it is a stationary solution
to the Cauchy problem for the TVF.

\noindent(iv) As a direct consequence of Theorem \ref{local.cont.modulus}, $C^{0,\alpha}$-regularity is preserved along the flow for any $\alpha\in (0,1]$. (Similar results have been obtained for the denoising problem and for the Neumann problem for the TVF in \cite{CCN}.)
Moreover, if $u_0 \in W^{1,1}(\RR)$, then $u(t) \in W^{1,1}(\R)$ (this is a consequence of the fact that the oscillation does not increase on any subinterval).

\noindent(v) If $u_0 \in BV_{loc}(\RR)$, a priori we do not have a well-defined semigroup.
However, in this case $u_0$ is locally bounded and the set of its discontinuity points is countable  (see \cite[Section 3.2]{AFP}),
and so in particular has Lebesgue measure zero.
Hence, by classical theorems on the Riemann integrability of functions,
we can find two sequences of step functions such that
$u_0^{\e,-} \leq u_0 \leq u_0^{\e,+}$ and $\|u_0^{\e,+}-u_0^{\e,-}\|_1\leq \e$
{(the number of steps will in general be infinite,
but finite on any bounded interval)}.
Then, by approximation we can still define a dynamics, which will still be
contractive in any $L^p$ space.

\noindent\textbf{Behaviour near maxima and minima. }
We conclude this section by giving an informal description of the evolution of a general solution
(excluding ``pathological'' cases).

Assume that $u_0$ has a local maximum at $x_0$.
Then, at least for short time, the solution is explicitly given near $x_0$ by
\[
u(t,x)=\min\{u_0(x),h(t)\},
\]
where the constant value $h(t)$ is implicitly defined by
\[
\int_{I_0} \big[ u_0(x) - h(t)\big]_+\dx= 2t\,,
\]
$I_0$ being the connected component of $\{u_0>h(t)\}$ containing $x_0$,
see Figure \ref{fig.norate.1}.
For a minimum point the argument is analogous.
The dynamics goes on in this way until a local minimum ``merges'' with a local maximum,
and then one can simply start again the above description starting from the new configuration.
\begin{figure}[ht]
\includegraphics[height=5cm, width=6cm]{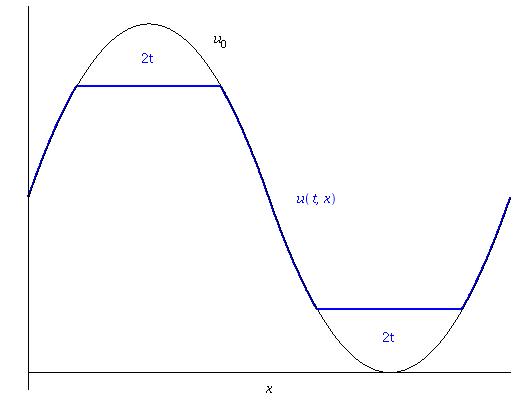}\includegraphics[height=5cm, width=9.75cm]{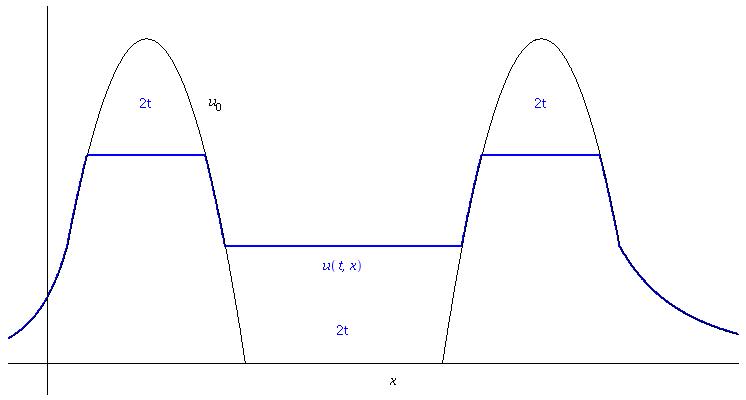}
\caption{\small \textit{Dynamic of TVF at maximum and minimum points.}}
\label{fig.minmax.cont}
\end{figure}

\subsection{Rescaled flow and stationary solutions}\label{resc.flow.sect}
Let $u_0 \in BV(\R)\cap L^1(\R)$ be a nonnegative compactly supported initial datum.
First we show that $u_0$ extinguishes in finite time, and we calculate the explicit
extinction time. (Note that, even in general dimension, estimates from above and from below on the extinction time were already known, see for example \cite{ACDM, ACM, Giga9}.)

\begin{proposition}[Loss of mass and extinction time]\label{ext.time+mass.loss}
Let $u(t)$ be the solution to the Cauchy problem in $\R$ for the TVF, starting from a non-negative compactly supported initial datum $u_0\in L^1(\RR)$. Then the following estimates hold:
\begin{equation}\label{loss.mass}
\int_\RR u(t,x)\dx=\int_\RR u_0(x)\dx -2t= 2(T-t)\qquad\mbox{for all $t\ge 0$}\,,
\end{equation}
and the extinction time for $u$ is given by
\begin{equation}\label{ext.time}
T=T(u_0)=\frac{1}{2}\int_\RR u_0(x)\dx\,.
\end{equation}
\end{proposition}
\noindent {\sl Proof.~}
Arguing by approximation and using the stability of the TVF in $L^1$, it suffices to consider the case when $u_0$
is a nonnegative step function. Assume that $\supp(u_0)\subset [a,b]$
and that $u_0$ jumps both at $a$ and at $b$. Then, by the explicit formula for $u(t)$,
we immediately deduce that $u(t)$ jumps both at $a$ and $b$, with $z(t,a)=1$
and $z(t,b)=-1$ (since $u(t)$ is nonnegative as well). Hence
\[
\frac{\rd}{\dt}\int_\RR u(t,x)\dx=\int_a^b \partial_x z(x,t)\dx= z(b)-z(a)=-2,
\]
from which \eqref{loss.mass}-\eqref{ext.time} follow.\qed

\noindent \textbf{Remark. }Let us point out that there is no general explicit formula for the extinction time when $u_0$ changes sign.

\noindent \textit{The rescaled flow. }We now are interested in describing the behavior of the solution near the extinction time.  To this end we need to perform a logarithmic time rescaling, which maps the interval $[0,T)$ into $[0,+\infty)$, where $T$ is the extinction time corresponding to the initial datum $u_0$. We define
\begin{equation}\label{TVF.rescaled}
w(s, x)=\frac{T}{T-t}\,u\left(t,x\right)\,,\qquad Z(s,x)=z(t,x)\,,\qquad s =T\log\left(\frac{T}{T-t}\right)\,,\qquad t=T\left(1-\ee^{-s/T}\right)\,,
\end{equation}
where $u(t)$ is a solution to the TVF. Then
\begin{equation}
\label{eq:rescaled}
\partial_s w(s,x)=\partial_x Z+\frac{w}{T},\qquad Z\cdot D_x w=|D_x w|\,,\qquad w(0,x)=u_0(x).
\end{equation}
We observe that stationary solutions $S(x)$ for the rescaled equation for $w$ correspond to separation of variable solutions in the original variable, namely
\[
-\partial_x Z=\frac{S}{T}\qquad\mbox{provides the separate variable solution}\qquad U_T(t,x):=\frac{T-t}{T}S(x)\,.
\]
We need now to characterize the stationary solutions. To this aim, we first have to define the ``extended support''of a function $f$
as the smallest interval that includes the support of $f$:
\[
\supp^*(f)=\inf\left\{[a,b]\;|\; \supp(f)\subseteq [a,b]\right\}.
\]

\begin{theorem}[Stationary solutions]\label{thm.stat.sols}
All compactly supported solutions of the equation
\begin{equation}\label{stationary.eq}
-\partial_x Z=\frac{S}{T}, \qquad Z\cdot D_x S=|D_x S|,
\end{equation}
are of the form
\begin{equation}\label{stationary.ab}
S(x)=\frac{2T}{b-a}\,\chi_{[a,b]}(x)\,,
\end{equation}
with $[a,b]\subseteq\RR$.
\end{theorem}
\noindent {\sl Proof.~}
Let us assume that $\supp^*{S}= [a,b]$. Since $S$ is nonnegative
we have $Z(a)=1$, $Z(b)=-1$.
We claim that $-1<Z<1$ on $(a,b)$.

Indeed, assume by contradiction that $Z(x_0)=1$ for some point $x_0 \in (a,b)$
(the case resp. $Z(x_0)=-1$ is completely analogous).
Then, using again that $S$ is nonnegative we obtain
$\partial_x Z=-\frac{S}{T} \leq 0$, which implies that
$Z\equiv 1$ on $[a,x_0]$.
Hence $S=-T\,\partial_x Z=0$ on $[a,x_0]$, which contradicts the definition of $\supp^*(S)$.

Thanks to the claim, since $Z\cdot D_x S=|D_x S|$ we easily deduce that $|D_x S|(a,b)=0$,
that is, $S$ is constant inside $(a,b)$.
To find the value of such a constant, we simply integrate the equation over $[a,b]$,
and we get
\[
\int_a^b \frac{S}{T}\dx=-\int_a^b \partial_x Z(x)\dx=Z(a)-Z(b)=2.
\]
and \eqref{stationary.ab} follows.\qed

\begin{corollary}[Separate variable solutions]
All compactly supported solutions of the TVF  obtained by separation of variables are of the form
\begin{equation}\label{sep.var.sols.orig.time}
U_T(t,x)=2\frac{T-t}{b-a}\chi_{(a,b)}(x)\,.
\end{equation}
where $T>0$ and $[a,b]\subseteq\RR$.
\end{corollary}

\begin{proposition}[Mass conservation for rescaled solutions]
Let $u(t)$ be the solution to the TVF corresponding to a nonnegative initial datum
$u_0\in BV(\R) \cap L^1(\RR)$. Let $w(s)$ be the corresponding rescaled solution, as in \eqref{TVF.rescaled}, then we have that
\begin{equation}\label{mass.conservation}
\int_\RR w(s,x)\dx=\int_\RR u_0(x)\dx \,.
\end{equation}
\end{proposition}
\noindent {\sl Proof.~}From \eqref{loss.mass},  \eqref{TVF.rescaled}, and
 the fact that the extinction time is given by $2T=\int_\RR u_0\dx$, we deduce that
\[
\int_\RR w(s,x)\dx=\frac{T}{T-t}\left[\int_\RR u_0(x)\dx -2t\right]=\ee^{s/T}\left[\int_\RR u_0(x)\dx -2T\left(1-\ee^{-s/T}\right)\right]=\int_\RR u_0(x)\dx\,.
\]
\qed

\begin{proposition}[Stationary solutions are asymptotic profiles]\label{limits.prop}
Let $w(s,x)$ be a solution to the rescaled TVF corresponding to a non-negative initial datum $u_0\in BV(\RR)\cap L^1(\RR)$. Then there exists a subsequence $s_n\to \infty$ such that $w(s_n,\cdot)\to S$ in $L^1(I)$ as $n\to \infty$ where $S$ is a stationary solution as in \eqref{stationary.ab}. Equivalently we have that there exists a sequence of times $t_n\to T$ as $n\to\infty$ such that
\[
\left\|\frac{u(t_n,\cdot)}{T-t_n}-\frac{S}{T}\right\|_{L^1}\xrightarrow[n\to\, \infty]{} 0\,.
\]
where $S$ is a stationary solution.
\end{proposition}
\noindent {\sl Proof.~}This is a well known result, see e.g. Theorem 4.3 of \cite{ABC} or Theorem 3 of \cite{ACDM} for the homogeneous Dirichlet problem on bounded domains. See also the book \cite{ACM}\,.\qed

\noindent\textbf{Remark. } From the above result we cannot directly deduce the correct extinction profile for the TVF, since there is not uniqueness of the stationary state, as Theorem \ref{thm.stat.sols} shows. Indeed a priori there may exists different subsequences such that the solution $w$ approaches two different stationary states along the two subsequences. We shall prove in the next section that such phenomenon does not occur.

\subsection{Asymptotics of the TVF}\label{asymptotic.sect}
Here we want to characterize the asymptotic (or extinction) profile for solutions to the TVF
 in function of the non-negative initial datum $u_0\in BV(\RR)\cap L^1(\RR)$.

\begin{theorem}[Extinction profile for solutions to the TVF]\label{Thm.Asymptotic}
Let $u(t,x)$ be a solution to the TVF corresponding to a non negative initial datum $u_0\in BV(\RR)$ with $\supp^*(u_0)=[a,b]$, and set
$$
T=\frac{1}{2}\int_a^b u_0(x)\dx\,.
$$
 Then $\supp(u(t))=[a,b]$ for all $t \in (0,T)$ and
\begin{equation}\label{asympt.TVF}
\left\|\frac{u(t,\cdot)}{T-t}-2\frac{\chi_{[a,b]}}{b-a}\right\|_{L^1([a,b])}\xrightarrow[t\to\, T]{} 0\,.
\end{equation}
\end{theorem}
\noindent\textbf{Remarks. }\noindent(i) The above theorem shows to important facts: firstly, the support of the solution becomes instantaneously the ``extended support'' of the initial datum, which is the support of the extinction profile.
Secondly, on $[a,b]=\supp^*(u_0)$ we consider the quotient $u(t,x)/U_T(t,x)$, where $U_T$ is the separate variable solution $U_T(t,x)=(T-t)S(x)$,
and $S(x)=2\frac{\chi_{[a,b]}}{b-a}$ (see \eqref{stationary.ab} and \eqref{sep.var.sols.orig.time}).
It is interesting to point out that $U_T$ is explicitly characterized in function of the extinction time
(i.e $\frac{1}{2}\int u_0$) and of extended support of the initial datum.
Then \eqref{asympt.TVF} can be rewritten as
\[
\left\|\frac{u(t,\cdot)}{U_T(t,\cdot)}-1\right\|_{L^1([a,b])}\xrightarrow[t\to\, T]{} 0\,.
\]
(In literature this result is usually called \textit{convergence in relative error}.) Equivalently, $L^1$-norm of the difference decays at least with the rate
\[
\big\|u(t,\cdot)-U_T(t,\cdot)\big\|_{\LL^1(\RR)}\le o(T-t)\,.
\]
In the next paragraph we will prove that the $o(1)$ appearing in the above rate cannot be
quantified/improved, so that  the convergence result of Theorem \ref{Thm.Asymptotic} is sharp.

\noindent(ii) The result of the above theorem can be restated in terms of the rescaled flow of Subsection \ref{resc.flow.sect}: if $w(s,x)$ is the rescaled solution corresponding to $w(0,x)=u_0$ (see \ref{TVF.rescaled}), then the relative error $w(s,x)/S\to 1$ as $s\to\infty$ in $L^1(\supp^*(u_0))$.

\medskip

\noindent {\sl Proof.~}The proof is divided into several steps.

\noindent$\bullet~$\textsc{Step 1. }\textit{Compactness estimates. }
Let $z(x,t)$ be associated to the $u(t,x)$ (as in the definition of strong solution, see Section \ref{sols.def.1d}). We claim that
\begin{equation}\label{cpt.est}
\|\partial_x z(x,t)\|_2\le \frac{2\|u_0\|_2}{t}\,.
\end{equation}
The proof of this fact is quite standard in semigroup theory, once one observes that $\partial_t u(t,x)=\partial_x z_u(x,t)$, Indeed, the homogeneity of the semigroup implies that
\begin{equation}\label{u.t.est}
\|\partial_t u(t,x)\|_p\le \frac{2\|u_0\|_p}{t}\qquad \forall \, p \in [1,\infty]\,.
\end{equation}
Although the latter estimate is classical and due to Benilan and Crandall \cite{BC2}, we briefly recall here the proof for convenience of the reader. Since $u(t,x)$ is a solution to the TVF starting from $u(0,\cdot)=u_0(\cdot)$\,, it is then clear that $u_\lambda(t,x)=\lambda u(\lambda^{-1}\,t,x)$ is again a solution to the TVF
starting from $u_\lambda(0,\cdot)=\lambda u_0(\cdot)$\,, for all $\lambda\ge 0$. Then
\[\begin{split}
u(t+h,x)-u(t,x)&=\frac{t+h}{t}u_\lambda(t,x)-u(t,x)=\lambda^{-1}u_\lambda(t,x)-u(t,x)\\
&=(\lambda^{-1}-1)u_\lambda(t,x) +(u_\lambda(t,x)-u(t,x))
\end{split}\]
where we have defined $\lambda:=t/(t+h)>0$. We now use the contraction property of the TVF in any $L^p-$space to conclude that
\[\begin{split}
\|u(t+h,x)-u(t,x)\|_p  &\le(\lambda^{-1}-1)\|u_\lambda(t,x)\|_p + \|u_\lambda(t,x)-u(t,x)\|_p\\
&\le |\lambda^{-1}-1|\|u_0\|_p +\|u_\lambda(0,x)-u(0,x)\|_p\\
&=|\lambda^{-1}-1|\|u_0\|_p +|\lambda-1|\|u_0\|_p = \frac{h}{t}+\frac{h}{t+h}\|u_0\|_p\le \frac{2h}{t}\|u_0\|_p.
\end{split}
\]
Letting $h\to 0^+$, \eqref{u.t.est} follows\,.

\noindent$\bullet~$\textsc{Step 2. }\textit{Stability up to $T^-$.}
Since $u_0 \in BV(\RR)$, as in Subsection \ref{sect:further}, Property (v),
we can find two sequences of step functions
$u_0^{\e,-} \leq u_0 \leq u_0^{\e,+}$ such that $\|u_0^{\e,+}-u_0^{\e,-}\|_1\leq \e$.
In particular, by the formula for the extinction time, we deduce that
$$
T-\e/2 \leq T(u_0^{\e,-}) \leq T(u_0^{\e,+})\leq T+\e/2\,.
$$
Moreover, up to replacing $u_0^{\e,+}$ with $\min\{u_0^{\e,+},\|u_0\|_\infty \chi_{\supp^*(u_0)} \}$,
we have that
$$
\supp^*(u_0^{\e,-}),\supp^*(u_0^{\e,+}) \to \supp^*(u_0)\qquad \text{as }\e \to 0.
$$
Since the evolution of step function is explicit, it is immediately checked that
$$
\supp(u^{\e,-}(t))=\supp^*(u_0^{\e,-}),\qquad \supp(u^{\e,+}(t))=\supp^*(u_0^{\e,+})
$$
for $t \in (0,T-\e/2]$ (indeed, for any subinterval $I\subset\subset\supp^*(u_0^{\e,-})$
where $u_0^{\e,-}$ vanishes, $u^{\e,-}(t)$ becomes instantaneously positive, see
Figure \ref{fig.minmax.cont}. By the parabolic maximum principle, this implies
$$
\supp^*(u_0^{\e,-}) \subset \supp(u(t)) \subset \supp^*(u_0^{\e,+}) \qquad \text{for $t \in (0,T-\e/2]$},
$$
which by the arbitrariness of $\e$ implies that
\begin{equation}
\label{eq:stability support}
 \supp(u(t))= [a,b]  \qquad\forall\, t\in (0,T).
\end{equation}

\noindent$\bullet~$\textsc{Step 3. }\textit{Convergence of $z$ and shape of $u(t)$
before the extinction time. }
By \eqref{eq:stability support} and the fact that $u(t)$ is nonnegative,
we deduce that $z(t,a)=1$ and $z(t,b)=-1$ for all $t \in (0,T)$.
Estimates \eqref{cpt.est} imply that $z(t,\cdot)$ is uniformly bounded in $C^{1/2}([a,b])$, hence it is compact in $C^0([a,b])$.
Moreover, by \eqref{eq:rescaled}
and Proposition \ref{limits.prop} there exists a sequence $s_k=T\log\left(\frac{T}{T-\e_k}\right) \to \infty$ (i.e. $\e_k \to 0^+$)
such that $\partial_s w(s_k,x) \to 0$. Hence, up to extracting a subsequence, we deduce that
$$
\lim_{k \to \infty} z(T-\e_k)=1-\int_a^x \frac{S(y)}{T}\,dy \qquad \text{uniformly on }[a,b],
$$
where $S=\frac{2T}{\b-\a}\,\chi_{[\a,\b]}(x)$ is a stationary solution, with $a \leq \a \leq \b \leq b$
(observe that, by Step 2, the support can only shrink).

\noindent$\bullet~$\textsc{Step 4. }\textit{Shape of $u(t)$
before the extinction time. }
By Step 3 we know that $z(T-\e_k)$ converges uniformly on $[a,b]$ to the function
\[
z_S(x)=\left\{
\begin{array}{ll}
1&\qquad\mbox{if }a\le x\le \alpha\\
\frac{\alpha+\beta}{\beta-\alpha}-\frac{2}{\beta-\alpha}x &\qquad\mbox{if }\alpha<x<\beta\\
-1&\qquad\mbox{if }\beta\le x\le b.\\
\end{array}
\right.
\]
Hence, for $k$ sufficiently large, there exist $\a<\a_{\e_k}<\b_{\e_k}<\b$
such that $-1<z(T-\e_k)<1$ on $[\a_{\e_k},\b_{\e_k}]$,
$z(T-\e_k)>-1$ on $[a,\a_{\e_k}]$, $z(T-\e_k)<1$ on $[\b_{\e_k},b]$,
and  $|\a-\a_{\e_k}|+|\b-\b_{\e_k}| \to 0$ as $\e \to 0$.
Since $z(T-\e_k)\cdot D_xu(T-\e_k)=|D_xu(T-\e_k)|$,
we easily deduce that
$u(T-\e)$ is increasing on $[a,\a_{\e_k}]$, constant on $[\a_{\e_k},\b_{\e_k}]$
and decreasing on $[\b_{\e_k},b]$.

\noindent$\bullet~$\textsc{Step 5. }\textit{Solutions with only one maximum point. }
Fix $k$ large enough, and consider the evolution of $u(t)$ on the time interval $[T-\e_k,T]$.
By Step 4 and the discussion at the end of Subsection \ref{sect:further}.
the evolution of $u(t)$ is explicit:
\[
u(t,x)=\min\{u(T-{\e_k},x),h(t)\}\qquad \forall \,t \in [T-\e_k,T)
\]
where $h(t)>0$ is implicitly defined by
\begin{equation}\label{loss.mass.2}
\int_\RR \big[ u(T-{\e_k},x) - h(t)\big]_+\dx= 2t\,.
\end{equation}
Since $\supp(u(t))=[a,b]$ (see Step 2),
the above formula shows that the set $[a(t),b(t)]$ where $u(t)$ is constant expands in time
and converges to $[a,b]$.
Moreover, from \eqref{loss.mass.2} we easily obtain the estimate
$$
-\frac{2}{b-a} \leq \dot h(t) \leq -\frac{1}{b-a}
$$
for $t$ close to $T$. Hence, the equation
$\partial_t u=\partial_x z$ implies that
\begin{equation}
\label{eq:strict decr}
-\frac{2}{b-a} \leq \partial_x z(t)\leq -\frac{1}{b-a}\qquad \text{on }[a(t),b(t)].
\end{equation}
Since $[a(t),b(t)] \to[a,b]$, the uniform convergence of $z(t)$ to $z_S$
is compatible with \eqref{eq:strict decr} if and only if
\[
z_S(x)=
\frac{a+b}{b-a}-\frac{2}{b-a}x\qquad\text{on }[a,b],
\]
i.e. the unique possible limiting profile is
$2\,T\frac{\chi_{[a,b]}}{b-a}$, as desired.
\mbox{\qed}

\subsection{Rates of convergence}\label{sect.rates}
After proving asymptotic convergence to a stationary state, the next natural question
is whether there exists a  universal rate of convergence to it.
As the next theorem shows,
the answer is negative.

\noindent Before stating the result, let us make precise what do we mean by decay rate.
\begin{definition} Let $\xi:[0,\infty)\to[0,\infty)$ be a continuous increasing function,
with $\xi(0)=0$. We say that $\xi$ is a \textit{rate function} if, for any solution $u(t)$ of the TVF,
\begin{equation}\label{no.rate.1}
\left\|\frac{\,u(t)}{T-t}-\frac{S}{T}\right\|_{L^1(I)}\le \xi(T-t)\qquad\mbox{for any $t$ close to the extinction time $T$\,. }
\end{equation}
\end{definition}
The following result shows that there cannot be a universal rate of convergence to any stationary profile.
\begin{theorem}[Absence of universal convergence rates]
For any rate function $\xi:[0,\infty)\to[0,\infty)$\,, there exists an initial datum $u_0\in BV(\RR)$,
with $\supp^*(u_0)=[0,1]$, such that
\begin{equation}\label{no.rate.1}
2\,\xi(T-t)\le \left\|\frac{u(t)}{T-t}-2\chi_{[0,1]}\right\|_{L^1(I)}\,,\qquad\mbox{for any }0\le T-t\le 1.
\end{equation}

\end{theorem}
\noindent {\sl Proof.~} Let us fix a rate function $\xi:\RR\to\RR$. It is not restrictive to assume that
$\xi$ is strictly increasing, and that $\xi(s)\ge s$ for any $s\in[0,1]$.

\begin{figure}[ht]
\centering
\includegraphics[height=5cm, width=7cm]{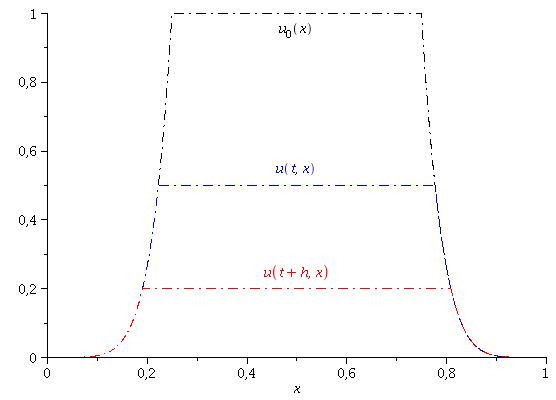}\includegraphics[height=5cm, width=8cm]{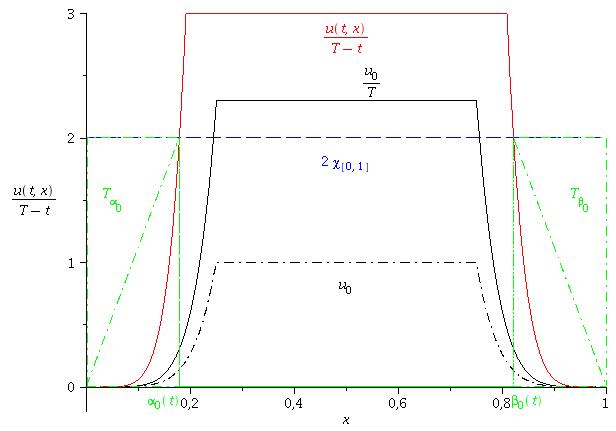}
\caption{\small \textbf{Left: }\textit{Dynamic of $u(t)$: } black: $u_0(x)$\,, \textcolor{darkblue}{blue: $u(t,x)$}\,, \textcolor{darkred}{red: $u(t+h)$}\,.
\textbf{ Right: }\textit{Rescaled dynamic: }\textrm{ black: $u_0(x)$ (dashdot) and $u_0/T$ (cont.), \textcolor{darkblue}{blue: $S(x)=2\chi_{[0,1]}$, } \textcolor{darkred}{red: $u(t,x)/(T-t)$, } \textcolor{darkgreen}{green: definition of $\alpha_0(t), \beta_0(t)$ and $\mathbb{T}_{\alpha_0}\,,\,\mathbb{T}_{\beta_0}$\,.}} }
\label{fig.norate.1}
\end{figure}

Let $\xi^{-1}=[0,\infty)\to[0,\infty)$ be the inverse of $\xi$, so that $\xi^{-1}(s)\le s$ for any $s\in[0,1]$, and choose the initial datum $u_0$ to be
\begin{equation}\label{u0.norate}
u_0(x)=\left\{
\begin{array}{lll}
c_0\xi^{-1}(x)\,,&\qquad\mbox{if }0\le x\le \frac{1}{4}\\
1\,,&\qquad\mbox{if }\frac{1}{4}<x<\frac{3}{4}\\
c_0\xi^{-1}(1-x)\,,&\qquad\mbox{if }\frac{3}{4}\le x\le 1\\
\end{array}
\right.
\end{equation}
with $c_0:=1/\xi^{-1}(1/4)$ (so that $u_0(1/4)=u_0(3/4)=1$).
By Theorem \ref{Thm.Asymptotic} we know that the solution $u(t)$ corresponding to $u_0$ extinguish at time
\[
\frac{1}{4}\le T=\frac{1}{2}\int_0^1u_0(x)\dx\le \frac{1}{2}
\]
and that $u(t)/(T-t)$ converges strongly in $L^1([0,1])$ to $S/T=2\chi_{[0,1]}$ as $t\to T$. First we prove that the $\LL^\infty$-norm satisfies the bound
\begin{equation}\label{infty.norm.bounds}
2(T-t)\le \|u(t)\|_\infty \le 4(T-t)\,.
\end{equation}
The first equality follows by the loss of mass formula \eqref{loss.mass} and H\"older inequality on
$[0,1]$:
\[
2(T-t)=\|u(t)\|_1\le \|u(t)\|_\infty\,.
\]
The second inequality follows as we know the explicit behaviour of the solution around maximum points (see the end of Subsection \ref{sect:further}), namely
\[
u(t,x)=\min\{u_0(x),h(t)\}
\]
where $h(t)>0$ satisfies
\begin{equation}\label{loss.mass.2.norate}
\int_0^1\big[ u_0(x) - h(t)\big]_+\dx= 2t\,.
\end{equation}
Then $u(t)$ is constant on an interval of the form $[\alpha(t),\beta(t)]$,
with $\alpha(t) \to 0^+$ and $\beta(t) \to 1^-$, see Figure \ref{fig.norate.1}.

\noindent Hence
\[
2(T-t)=\|u(t)\|_1=\int_0^1 u(t,x)\dx\ge \int_{1/2}^{3/4} u(t,x)\dx=\frac{1}{2}\|u(t)\|_\infty\,,
\]
and so $\|u(t)\|_\infty\le 4(T-t)$\,, and inequality \eqref{infty.norm.bounds} is proved.

Now, let $\alpha_0(t )\in[0,\alpha(t)]$ and $\b_0(t )\in[\b(t),1]$ be the unique points such that $u(t ,\alpha_0(t ))/(T-t )=u(t ,\b_0(t ))/(T-t )=2$
(such points exists thanks to the lower bound $\|u(t)\|_\infty\ge 2(T-t)$, see also Figure \ref{fig.norate.1}).
Let $\mathbb{T}_{\alpha_0}$ (resp. $\mathbb{T}_{\beta_0}$)
be the rectangular triangle with height $\{0\}\times [0,2]$ (resp. $\{1\}\times [0,2]$)
and basis $[0,\a(t)]\times \{2\}$ (resp. $[\b(t),1]\times \{2\}$),
as depicted in Figure \ref{fig.norate.1}. Denote by $|\mathbb{T}_{\alpha_0}|$,
$|\mathbb{T}_{\beta_0}|$ their measure.
Then, since $\xi^{-1}(s)\geq s$ we easily obtain the estimate
\[
\begin{split}
&\left\|\frac{u(t )}{T-t }-\frac{S}{T}\right\|_{L^1([0,1])}
=\int_0^1\;\left|\frac{u(t ,x)}{T-t }-2\chi_{[0,1]}(x)\right|\;\dx\\
&=\int_0^{\alpha(t )}\,\left|\frac{u(t ,x)}{T-t }-2\right|\,\dx
 +\int_{\alpha(t )}^{\beta(t )}\;\left|\frac{u(t ,x)}{T-t }-2\right|\,\dx
 +\int_{\beta(t )}^1\,\left|\frac{u(t ,x)}{T-t }-2\right|\;\dx\\
&\ge\int_0^{\alpha(t )}\,\left|\frac{u(t ,x)}{T-t }-2\right|\,\dx
 +\int_{\beta(t )}^1\,\left|\frac{u(t ,x)}{T-t }-2\right|\,\dx\\
 &\ge\int_0^{\alpha_0(t )}\,\left(2-\frac{u(t ,x)}{T-t }\right)\,\dx
 +\int_{\beta_0(t )}^1\,\left(2-\frac{u(t ,x)}{T-t }\right)\,\dx\\
&\ge |\mathbb{T}_{\alpha_0} |+ |\mathbb{T}_{\beta_0}|=2\alpha_0(t ).\\
\end{split}
\]
To estimate $\a_0(t)$ from below, we observe that on $[0,\alpha(t)]$ we have that  $u(t ,\alpha_0(t ))=u_0(\alpha_0(t ))=c_0\xi^{-1}(\alpha_0(t ))$, so
\[
\alpha_0(t )=\xi\left(\frac{2}{c_0}(T-t )\right).
\]
Now, recalling that $\xi$ is strictly increasing and $\xi(1)\geq 1$, we get $\xi(2)\geq 1/4$,
or equivalently $2/c_0=2/\xi^{-1}(1/4)\ge 1$.
Hence $\alpha_0(t) \geq \xi(T-t)$, which concludes the proof.
\qed

\noindent\textbf{Remark. }The above Theorem shows that there cannot be universal rates of convergence. A similar construction will provide (nontrivial) initial data
for which the convergence is as fast as desired.

\begin{theorem}[Fast decaying initial data]
For any rate function $\xi:[0,\infty)\to[0,\infty)$\,, there exists an initial datum $u_0\in L^1(I)$ such that the corresponding solution $u(t)$ satisfies
\begin{equation}\label{no.rate.2}
\left\|\frac{u(t)}{T-t}-2\chi_{[0,1]}\right\|_{L^1(I)}\le \,\xi\big(8(T-t)\big)\,,\qquad\mbox{for any }0\le T-t\le 1\,.
\end{equation}
\end{theorem}
\noindent {\sl Proof.~} Fix a rate function $\xi:\RR\to\RR$, which is continuous, increasing, $\xi(0)=0$, and $\xi(s)\le s$.

\begin{figure}[ht]
\centering
\includegraphics[height=5cm, width=7cm]{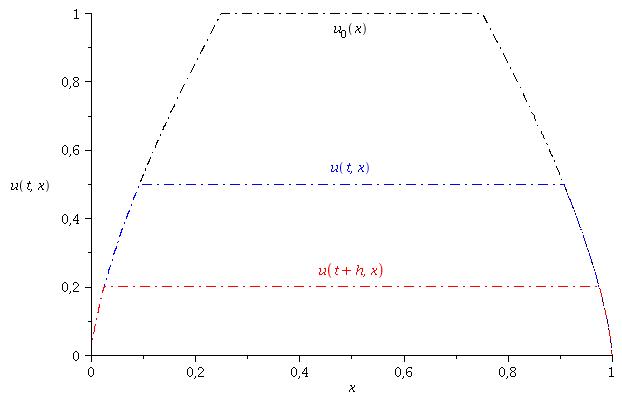}\includegraphics[height=5cm, width=9cm]{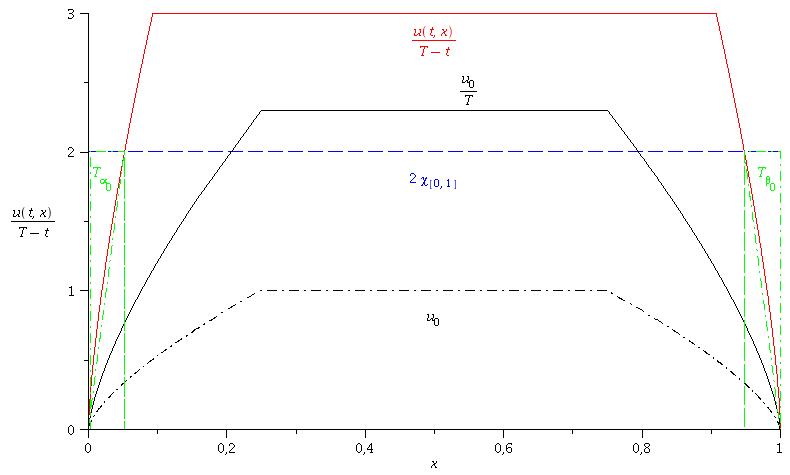}
\caption{\small \textbf{Left: }\textit{Dynamic of $u(t)$: } black: $u_0(x)$\,, {blue: $u(t,x)$}\,, \textcolor{darkred}{red: $u(t+h,x)$}\,.
\textbf{ Right: }\textit{Rescaled dynamic: }\textrm{ black: $u_0(x)$, {blue: $S(x)=2\chi_{[0,1]}$, } \textcolor{darkred}{red: $u(t,x)/(T-t)$, } \textcolor{darkgreen}{green: definition of $\alpha_0(t), \beta_0(t)$ and $\mathbb{T}_{\alpha_0}\,,\,\mathbb{T}_{\beta_0}$\,.}} }
\label{fig.norate.3}
\end{figure}

Let $\xi^{-1}=[0,\infty)\to[0,\infty)$ denote its inverse,
and define $u_0$ as in \eqref{u0.norate}, see Figure \ref{fig.norate.3}.
Then an analysis analogous to the one done in the previous Theorem proves the result.
We leave the details to the interested reader.\qed

\section{Solutions to the SFDE and solutions to the TVF}
\label{sect:TVF fast diff}

As explained in the introduction, TVF and SFDE are formally related by the fact that
``$u$ solves the TVF if and only if $D_xu$ solves the SFDE''.

In order to make this rigorous, we need first to explain what do we mean by a solution
of the SFDE, and then we will prove the above relation by approximating the TVF with the $p$-Laplacian and the SFDE by the porous medium equation.

The notion of solution we consider for the SFDE is the one of mild solution.
More precisely, since the multivalued graph of the function $r \mapsto \sign(r)$
is maximal monotone, by the results of Benilan and Crandall \cite{BC}, there exists
a continuous semigroup $S_0^t:L^1(\RR)\to L^1(\RR)$ such that
$S_t^0 v_0\in C([0,\infty);L^1(\R))$ is a mild solution of the SFDE. To be more precise, let $\varphi$ be a maximal monotone graph in $\RR$ (see \cite{B1}) and consider the problem
\begin{equation}\label{probl.SFDE}
\left\{\begin{array}{lll}
u_t=\Delta\varphi(u)\,, &\qquad\mbox{in }\mathcal{D}'\big((0,\infty)\times\RR\big)\\
u(0,x)=u_0(x)\,, & \qquad x\in \RR\\
\end{array}
\right.
\end{equation}
where the first equation is meant in the sense that
\begin{equation}\label{mild.sol.SFDE}
u_t=\Delta w \qquad\mbox{in }\mathcal{D}'\big((0,\infty)\times\RR\big)\,,\qquad\mbox{with }w(t,x)\in \varphi(u(t,x))\qquad\mbox{a.e. }t, x\in \RR\,.
\end{equation}
We now recall the celebrated results of Benilan and Crandall \cite{BC} adapted to our setting, namely it is sufficient to consider $\varphi(r)=\sign(r)$.
\begin{theorem}[Benilan-Crandall, \cite{BC}]\label{BC-THM}
Under the running assumptions, the following results hold true:\\
(i) There exists a unique solution $u\in C([0,\infty);L^1(\RR))\cap L^\infty([0,\infty)\times\RR)$ corresponding to the initial datum $u_0\in L^1(\RR)\cap L^\infty(\RR)$ such that \eqref{probl.SFDE} and \eqref{mild.sol.SFDE} are satisfied.\\
(ii) Let $u_n\in C([0,\infty);L^1(\RR))\cap L^\infty([0,\infty)\times\RR)$ be solutions of \eqref{mild.sol.SFDE} corresponding to the sequence $\varphi_n:\R \to \R$, $n=1,2,\ldots$ of maximal monotone graphs in $\RR$. Assume that $0\in\varphi_n(0)$,
\begin{equation}\label{hip.BC}
\lim_{n\to \infty}\varphi_n(r)=\varphi(r)\quad\mbox{for }r\in\RR,\qquad\mbox{and}\qquad \lim_{n\to \infty}\|u_{0,n}-u_0\|_{\LL^1(\RR)}=0\,.
\end{equation}
Then $u_n\to u$ in $C\big([0,\infty);L^1(\RR)\big)$\,, where $u$ is the solutions of \eqref{mild.sol.SFDE} corresponding to $\varphi$.
\end{theorem}
Now, let $S^m_t$ be the semigroup associated to the FDE equation
$$
\partial_t v=\Delta (v^m).
$$
Since the graphs of the function $r \mapsto r^m:=|r|^{m-1}r$ converge to the graph of $r \mapsto \sign(r)$, we can use Theorem \ref{BC-THM} with the simple choice $\varphi_n(r)=|r|^{\frac{1}{n}-1}r$ to guarantee that we have convergence (indeed as $m \to 0^+$) of $S_t^m v_0$ to $S_t^0 v_0$ in $C([0,\infty);L^1(\R))$, for any initial datum $v_0 \in L^1$.

On the other hand, we can consider the $p$-Laplacian semigroup $T_t^p$ for $p=1+m$.
It is well known {(see for example \cite{B1} chap. 4, or \cite{ACM} chap. 5)} that if $u_0 \in W^{1,p} (\RR)$,
then $T_t^p u_0 \in W^{1,p}(\RR)$, so that as $p \to 1^+$, strong solutions to the $p$-Laplacian converge to strong solutions to the TVF. (A detailed proof of this fact in dimension $n\ge 1$ can be found for instance in \cite{ACDM}, pg. 138-142, in the framework of the Dirichlet problems on bounded domains, but it can be easily adapted to other problems, including the Cauchy one.) Hence, if $u_0$ is a smooth compactly supported function, $T_t^p u_0 \to T_t^1 u_0$ in $C\big([0,\infty);L^1(\RR)\big)$ as $p \to 1^+$, where $T_t^1$ denotes the TVF-semigroup.

Moreover, if $p=1+m$, we have that $\partial_x \bigl( T_t^p u_0\bigr)$ solves (in the distributional and semigroup sense) the FDE with initial datum $\partial_x u_0$, i.e. $\partial_x \bigl( T_t^p u_0\bigr)=S_t^m \bigl( \partial_x u_0\bigr)$.
Hence, by letting $m \to 0^+$, we recover such a relation in the limit $p=1$ and $m=0$
by what we said above.
We can summarize this discussion in the following:
\begin{theorem}
Assume $u_0$ is a smooth compactly supported function. Let $1<p\le 2$, and $m=p-1$.
Then the following diagram is commutative:
\begin{diagram}[labelstyle=\textstyle]
T^p_t u_0\in W^{1,p}(\RR)\; & \rTo_{}^{\quad p\to 1^+ \quad } & \; T^1_t u_0\in W^{1,1}(\RR) \\
\dTo<{\partial_x} &   & \dTo>{\partial_x} \\
S^m_t \bigl(\partial_x u_0\bigr)\in L^{1+m}(\RR)\;  & \rTo_{\quad m\to 0^+\quad} & \; S^0_t \bigl(\partial_x u_0\bigr)\in L^1(\RR).
\end{diagram}
Note that the convergence in meant in the sense of distributions
\end{theorem}

At this point it is worth noticing that the vector field $z$ associated to the solution $u$ of the TVF as in Section \ref{sect.sol.TVF} and the function $w$ associated to the solutions to the SFDE as in \eqref{mild.sol.SFDE} are the same (just by letting $v_0= \partial_x u_0$ and $v=\partial_x u$)
\[
\partial_t u= \partial_x z \qquad \xrightarrow[~~~~\partial_x~~~~]{}\qquad \partial_t v =\partial_t\partial_x u=\partial_{xx} z= \partial_{xx}w\,.
\]

\noindent\textbf{Measures as initial data.}  Once the correspondence between TVF and SFDE is established for smooth initial data,
by stability in $L^1$ of both semigroups
it immediately extends to $u_0 \in W^{1,1}(\RR)$,
and then by approximation to $BV(\RR)\cap L^1(\R)$
initial data. However, at the level of the SFDE this would correspond to finite measures
$v_0$ such that $\int_{-\infty}^x v_0(dy) \in L^1(\RR)$, which is possible if and only if
$\int_{-\infty}^{+\infty} v_0(dy)=0$.
Actually, this class of data correspond exactly to the one for which there is extinction in finite time
(as this is the case for $L^1$ initial data to the TVF).

To remove this unnatural constraint on $v_0$, we observe that, by Subsection \ref{sect:further}, Property (v),
the TVF defines a contractive semigroup also on initial data which are only in $BV_{loc}(\RR)$.
In particular, the TVF is well-defined on data of the form $u_0(x)=\int_{-\infty}^x v_0(dy)$,
where $v_0$ is a (locally) finite measure on $\R$. Hence, this allows to define measure valued solutions of the SFDE as $\partial_x T_t^1(u_0)$, and this notion coincides with the one that one would get by considering
weak$^*$ limit of $L^1$ solutions.

\noindent Summing up, we have shown that:

\noindent$\bullet$ If $v_0 \in L^1(\RR)$, the unique mild solution { of the SFDE of Theorem \ref{BC-THM} is given by}
\begin{equation}
\label{eq:TVF SFDE}
S_t^0v_0=\partial_x \biggl(T_t^1\Bigl(\int_{-\infty}^x v_0(dy)\Bigr) \biggr).
\end{equation}

\noindent$\bullet$ Using \eqref{eq:TVF SFDE} we can uniquely extend the generator $S_t^0$
to measure initial data (actually, since the semigroup $T_t^1$ is well-defined on $L^2(\RR)$,
one could even extend the SFDE to distributional initial data in $W^{-1,2}(\RR)$).

\subsection{The 1-dimensional SFDE}\label{SFDE.sect}

By \eqref{eq:TVF SFDE}, the dynamics of the SFDE can be completely recovered from the one of the TVF.
We begin by illustrating some basic properties, and then, instead of trying to give a complete description of the evolution (which, by the analysis of the TVF done in the previous sections, would just be a tedious exercise), we prefer to briefly illustrate
the evolution of solutions of the SFDE in some simple but representative situations. As pointed out in the previous section, we are allowed to consider measures as initial data for the SFDE, keeping in mind that the distributional $x$-derivative of a solution to the TVF is a solution to the SFDE.

In the same way as step functions allowed us to understand the dynamics of the TVF,
we start by considering the dynamics of sum of delta masses, which is directly deduced from the one of step functions for the TVF.

\noindent\textbf{Example 1. Delta masses as initial data. }
Let us assume that $v_0=\sum_{i=1}^Na_i \delta_{x_i}$, with $x_1<\dots<x_N$. Then, for $t>0$ small (the smallness depending on the size of $|a_i|$)
\[
v(t)=\sum_{i=1}^Na_i(t) \delta_{x_i},
\]
with $a_i(0)=a_i$ and
\begin{equation}\label{ai}
a_i(t)=\left\{
\begin{array}{lll}
a_i\,,&\qquad\mbox{if }\sign(a_{i-1})=\sign(a_{i+1})=\sign(a_i)\\
\sign(a_i)\big(|a_i|-4t\big)_+\,,&\qquad\mbox{if }\sign(a_{i-1})=\sign(a_{i+1})=-\sign(a_i)\\
\sign(a_i)\big(|a_i|-2t\big)_+\,,&\qquad\mbox{if }\sign(a_{i-1})\,\sign(a_{i+1})=-1,\\
\end{array}
\right.
\end{equation}
where we use the convention $\sign(a_0):=\sign(a_1)$ and $\sign(a_{N+1}):=\sign(a_N)$.
This formula holds true until one mass disappear at some time $t_1'>0$, and then it suffices to
$v(t'_1)$ as initial data and repeat the construction (compare with Subsection \ref{dyn.step.funct}).

We observe that if all $a_i$ have the same sign, then $v_0$ produces a stationary solution $v(t,x)=v_0(x)$.
Moreover, the total mass is conserved under the dynamics.
In particular, $v(t)$ extinguishes in finite time if and only it has zero mean, i.e. $\sum_i a_i=0$
(however, there is no simple formula for the extinction time, see the remark after Proposition
\ref{ext.time+mass.loss}).\\

\noindent\textbf{General properties. }Arguing by approximation (or again using the direct relation with the TVF), as a consequence we have the following properties of the SFDE flow:

\noindent$(i)$ \textit{Nonnegative initial data. }Let  $v_0\ge 0$ be a locally finite measure, and define $u_0(x):=\int_{-\infty}^x v_0(dy)\ge 0$. Since $u_0$ is monotone non-decreasing,
it does not evolve under the TVF, cf. (iii) in Section \ref{sect:further}.
Hence $v_0$ is a stationary solution to the SFDE.
(Actually, since monotone profiles are the only stationary state for the TVF, the only stationary solutions for the SFDE are nonnegative/nonpositive initial data.)

\noindent$(ii)$ \textit{Only zero mean valued initial data extinguish in finite time. }
If $v_0$ is a finite measure, $v(t)$ converges in finite time to a stationary solution $\bar v$
such that $\int_\R \bar v(dy)=\int_\R v_0(dy)$. Moreover, $\bar v\equiv 0$ (i.e. $v_0$ extinguish in finite time) if and only if  $\int_\R v_0(dy)=0$.\\

\noindent\textbf{Example 2. Interaction between a delta and a continuous part. }

\begin{figure}
\centering
\includegraphics[height=7cm, width=14cm]{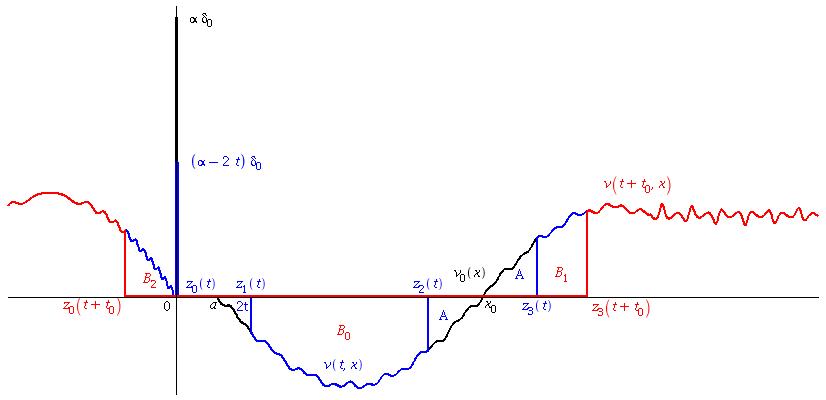}
\caption{black: $v_0(x)$\,, \textcolor{darkblue}{blue: $v(t,x)$, $t<t_0$}, \textcolor{darkred}{red: $v(t_1,x)$ is a stationary state}.}
\label{fig.SFDE.2}
\end{figure}

Let $v_0=\hat{v}_0+\alpha\delta_0$ where $\hat{v}_0\in C(\RR)$ is positive exception made for an interval $[0,x_0]$ as depicted in Figure \ref{fig.SFDE.2}. The zero set of $\hat{v}_0$ is $Z(\hat{v}_0)=[0,a]\cup\{x_0\}$, and a Dirac's delta with mass $\alpha$ is put at $x=0$.
We assume $\alpha$ is much smaller than $\int_a^{x_0} |\hat v_0|$.
The evolution basically changes at two steps:

\noindent$\bullet~$ $0\le t\le t_0$. The delta mass starts to lose its mass $\alpha$ by a factor $2t$ until time $t_0=\alpha/2$ (when it extinguish). This mass is compensated by
a ``gain of mass'' of $\hat v_0$: its zero set near $a$ starts to move to the right, and at time $t$ it is at a position $z_1(t)$ characterized by the fact that $\int_a^{a+z_1(t)}  |\hat v_0|=2t$ (cf. the blue area $2t$ in the Figure \ref{fig.SFDE.2}).
On the other hand, the isolated point $x_0\in Z(v_0)$ starts to ``expand'', creating a zero set
$[z_2(t),z_3(t)]\in Z(u(t))$, with $\int_{z_2(t)}^{x_0} |\hat v_0|=\int_{x_0}^{z_3(t)} |\hat v_0|=2t$.
This expansion of the zero set corresponds to the creation of flat parts at the level TVF (see Figure \ref{fig.minmax.cont}).

\noindent$\bullet~$ $t_0\leq t\leq t_1$. At time $t_0$ the delta disappears, and $u(t_0,x)$ is a piecewise continuous function on $\RR$ which is positive outside $[0,z_3(t_0)]$, it is zero on $Z(u(t_0,x))=[0,z_1(t_0)]\cup[z_2(t_0),z_3(t_0)]$, and is negative on $(z_1(t_0),z_2(t_0))$.
Starting form this time,
the free boundary expands on both components in the same way as described
above, until some time $t_1>t_0$ when $z_1(t_1)=z_2(t_1)$.
Observe that the loss of mass is such that $B_0=B_1+B_2$ (see Figure \ref{fig.SFDE.2}).

\noindent$\bullet~$ \textit{Reaching the stationary state in finite time. }
Since $u(t_1)\geq 0$, the solution becomes stationary
and $u(t)=u(t_1)$ for $t \geq t_1$.

\subsection{SFDE vs LFDE.}\label{sect:logFDE}
Let us go back to the fast diffusion equation \eqref{eq:FDE},
and assume that $v_0 \geq 0$ (so $v(t) \geq 0$ for all $t\geq 0$).
We remark that, changing the time scale $t \mapsto mt$,
the above equation can be written in two different ways which
lead to two different limiting equations: more precisely, setting $\rho(t,x)=v(t/m,x)$,
$$
\begin{array}{lll}
\partial_t v=\Delta\big(v^m\big) & \qquad\xrightarrow[m\to\, 0^+]\;\qquad & \partial_t v=\Delta\big(\sign(v)\big)\\[3mm]
\partial_t \rho =\div\big(\rho^{m-1}\nabla \rho\big) & \qquad\xrightarrow[m\to\, 0^+]\;\qquad & \partial_t \rho=\div\left(\rho^{-1}\nabla\rho\right)=\Delta \big(\log(\rho)\big)\,.\\
\end{array}
$$
Observe that
the evolution of $\rho(t,x)$ on the time interval $[0,T]$ corresponds to the evolution of $v(t,x)$
on the larger time interval $[0,T/m]$, and the diffusion
of $v$ is slower than the diffusion of $\rho$ by a factor $1/m$. So, when analyzing the limit as $m\to 0^+$, one gets two different limits, and the evolution of
$\rho(t,x)$ on the time interval $0\le t\le T$ corresponds to the evolution of $v(t,x)$ on the time interval
$0\le t<\infty$ for every $T>0$. This means that the solution to the SFDE corresponds to an evolution ``infinitely slower'' than the solution to the LFDE.

We saw that the Cauchy problem for the SFDE gives rise to a trivial dynamics on non-negative initial data.
However, the problem becomes non-trivial if we consider for instance the Dirichlet problem for the SFDE on a closed interval $I$ with zero boundary conditions. Indeed, by ``integrating in space'' such a solution, we obtain a solution to the TVF on $I$ with Neumann boundary conditions, whose dynamics on step functions has been described in Remark \ref{rmk:general bdry}.
In particular, from the example given there, we can explicitly find the dynamics of a finite
sum of positive deltas: if $v_0=\sum_{i=1}^N a_i \delta_{x_i}$ with $a_i >0$ and $x_1 \leq \ldots\leq x_N$, then $v(t)=v_0 - t[\delta_{x_1}+\delta_{x_N}]$ until one delta disappears, and then one simply restarts from there.
By approximation, we see that positive initial data
extinguish in finite time, and the extinction
time is given by $\frac{1}{2}\int_\RR v_0(dy)$.

This fact and the above discussion suggest that solutions to the LFDE should extinguish instantaneously, i.e. $\rho(t)\equiv 0$ for any $t >0$.
As one can deduce by the results of Rodriguez and V\'azquez \cite{RV2}, this is actually the case. (See also \cite{ERV, RV,RV2} and the book \cite{VazLN} for a complete theory of the logarithmic diffusion equation in one space dimension for positive initial data.)
So the above discussion is correct, and gives an heuristic explanation for this immediate extinction phenomenon.

\end{document}